\documentclass{amsart}
\usepackage{amssymb,amsthm}
\usepackage{pstricks}

\newcommand{\nc}{\newcommand}
\nc{\rnc}{\renewcommand}

\nc{\al}{\alpha}
\nc{\be}{\beta}
\nc{\ga}{\gamma}
\nc{\om}{\omega}
\nc{\Om}{\Omega}

\nc{\<}{\langle}
\rnc{\>}{\rangle}

\nc{\C}{\mathbf C}
\nc{\Z}{\mathbf Z}

\rnc{\P}{\mathbb P}

\nc{\inv}{^{-1}}
\nc{\xto}{\xrightarrow}
\nc{\SL}{\operatorname{SL}}
\rnc{\Sp}{\operatorname{Sp}}
\nc{\Norm}{\operatorname{N}}

\theoremstyle{definition}
\newtheorem{df}{Definition}[section]
\newtheorem{rem}[df]{Remark}
\newtheorem{ex}[df]{Example}

\theoremstyle{plain}
\newtheorem{theo}[df]{Theorem}
\newtheorem{prop}[df]{Proposition}

\newtheorem{cor}[df]{Corollary}

\numberwithin{equation}{section}

\rnc{\subjclassname}{%
  \textup{2000} Mathematics Subject Classification}

\title{A class of $T$-stable $(\P^1\times\dots\times\P^1)$'s in $G/B$}
\author{Christian Ohn}
\date{June 30, 2000; this version January 12, 2001}
\address{Universit\'e de Reims
\\
D\'epartement de Math\'ematiques (UPRESA 6056 du CNRS)
\\
Moulin de la Housse, B.P. 1039
\\
F-51687 Reims Cedex 2
\\
France}
\email{christian.ohn@univ-reims.fr}

\begin{document}

\begin{abstract}
Let $G$ be a connected complex semi-simple group, $B\subset G$ a Borel
subgroup, and $T\subset B$ a maximal torus. We construct a class of smooth
$T$-stable subvarieties inside the flag variety $G/B$, each of which is an
embedding of a product of projective lines.
\end{abstract}

\maketitle

\section{Introduction}

Let $G$ be a connected complex semi-simple group. Let $B$ be a Borel subgroup
and consider the natural (left) action of a maximal torus $T\subset B$ on the
flag variety $G/B$. A number of authors have studied $T$-stable subvarieties
in $G/B$; see e.g.~\cite{Br,Car,CarKur,CarKut,Dab,FlHa,GelSer,Kly}. In
particular, all $T$-stable curves are known~\cite{Car}: for each $w$ in the
Weyl group $W:=\Norm_G(T)/T$ and each root $\al$, there is a unique $T$-stable
curve through the $T$-fixed points $wB$ and $ws_\al B$. Each of these curves
is isomorphic to the projective line $\P^1$.

In this note, we describe a class of higher dimensional smooth $T$-stable
subvarieties in $G/B$, generalizing those $T$-stable curves. More precisely,
to each $w\in W$ and each set $\{\al_1,\dots,\al_d\}$ of pairwise orthogonal
roots (in the weaker sense, i.e.\ the sum of two $\al_k$'s may be a root), we
associate a $T$-stable subvariety in $G/B$ passing through all $T$-fixed
points $w'B$ with $w'$ of the form: $w$ times a product of some of the
$s_{\al_k}$. We then show (Theorem~\ref{maintheo}) that such a subvariety is
a closed embedding into $G/B$ of a product $\P^1\times\dots\times\P^1$ of $d$
projective lines. To the best of my knowledge, these varieties have not yet
been described in the literature (except for the curves mentioned above).
Also, I do not know whether they exhaust all $T$-stable subvarieties in $G/B$
that are isomorphic to a product of projective lines.

Although the varieties considered here are reminiscent of Bott-Samelson
varieties (see e.g.~\cite{BoSa,De,Mag}), we would like to point out some
differences: the latter are associated to certain sequences of simple roots
(not necessarily orthogonal), whereas our roots need not be simple (but must
be orthogonal). Our varieties are direct products of $\P^1$'s, whereas
Bott-Samelson varieties are (generally nontrivial) $\P^1$-fibrations over
Bott-Samelson varieties of lower dimension. (However, if the roots involved
are both simple \emph{and} pairwise orthogonal, then the successive
$\P^1$-fibrations become trivial and both constructions agree: we obtain the
Schubert variety corresponding to the product of the reflections associated
to these roots.)

Finally, let us remark that we found the present class of subvarieties of
$G/B$ while investigating an approach to quantum analogues of flag
varieties~\cite{qflag}. Since they could be of independent interest to
algebraic geometers, we decided to describe them in this note, separately
from~\cite{qflag}.

{\bf Acknowledgements.} I would like to thank M.~Brion and H.~Kraft for
useful discussions, as well as the Universit\'e de Reims for granting a
sabbatical leave during the year 1999--2000, when the present work has been
done.

\section{Notation}

\begin{tabular}{cl}
$G$&a connected complex semi-simple group
\\
$B,B^{-}$&two opposite Borel subgroups of $G$
\\
$T$&the maximal torus $B\cap B^{-}$
\\
$U,U^{-}$&the unipotent radicals of $B,B^{-}$
\\
$W$&the Weyl group $\Norm_G(T)/T$
\\
$\Phi$&the root system of $G$ w.r.t.\ $T$
\\
$\Phi^{+},\Phi^{-}$&the sets of positive and of negative roots w.r.t.\ $B$
\\
$<$&the partial order on $\Phi$ defined by $\al<\be\iff\be-\al\in\Phi^{+}$
\\
$s_\al$&the reflection in $W$ associated to a root $\al$
\\
$L_\al$&the copy of $(\mathrm P)\!\SL(2)$ in $G$ corresponding to a root
$\al$
\\
$U_\al$&the root group corresponding to a root $\al$
\\
$B_\al$&the Borel subgroup of $L_\al$ containing $U_\al$
\end{tabular}

\section{A class of subvarieties in $G/B$}

\begin{df}
An \emph{orthocell} in $W$ is a left coset in $W$ of a subgroup generated by
pairwise commuting reflections.
\end{df}
To each orthocell $C\subset W$, we will associate a $T$-stable subvariety
$E(C)\subset G/B$. To define it, we will have to make some choices: first,
choose $w\in C$, and choose a representative $\dot w\in\Norm_G(T)$ of $w$.
By definition, we have $C=w\<s_\al\mid\al\in\Om\>$ for some set $\Om$ of
positive and pairwise orthogonal roots: next, choose a numbering
$\al_1,\dots,\al_d$ of the elements of $\Om$ that is \emph{nonincreasing}, in
the sense that
\[
\al_k\not<\al_{k'}\qquad\text{for all $k<k'$.}
\]
(This is always possible: choose a maximal element in $\Om$ and call it
$\al_1$, then choose a maximal element among the remaining ones and call it
$\al_2$, etc.) Note that the numbering may be chosen arbitrarily unless
$\Phi$ has a component of type B$_n$, C$_n$, or F$_4$.

We then define
\[
E(C):=\{\dot wg_1\dots g_dB\mid g_k\in L_{\al_k}\,\forall k\}\subset G/B.
\]
In due course, we will show that $E(C)$ only depends on the orthocell $C$,
and not on the choices we have made above (see Remarks~\ref{indepnumb}
and~\ref{indepw}).
\begin{rem}
The set $E(C)$ contains all $T$-fixed points $w'B$ with $w'\in C$.
\end{rem}
\begin{prop}
\label{prodP1}
If $\al_1,\dots,\al_d$ are as above, then the map
\[
L_{\al_1}\times\dots\times L_{\al_d}\to G/B:(g_1,\dots,g_d)\mapsto \dot
wg_1\dots g_dB
\]
factors down to a map
\[
j:L_{\al_1}/B_{\al_1}\times\dots\times L_{\al_d}/B_{\al_d}\to G/B.
\]
\end{prop}
\begin{proof}
For all $1\le k\le d$, let $g_k\in L_{\al_k}$ and $b_k\in B_{\al_k}$. We need
to show that $g_1b_1\dots g_db_dB=g_1\dots g_dB$. By induction over $d$, we
may assume that the left hand side is equal to $g_1b_1g_2\dots g_dB$, so we
must show that $g_d\inv\dots g_2\inv b_1g_2\dots g_d\in B$.

We decompose $b_1=tu$, where $t\in B_{\al_1}\cap T$ and $u\in U_{\al_1}$.
Orthogonality of $\al_1,\dots,\al_d$ already implies that for each $2\le k\le
d$, $t$ commutes with all elements of $U_{\al_k}$ and of $U_{-\al_k}$. Since
$L_{\al_k}=\<U_{\al_k},U_{-\al_k}\>$ for each $k$, it follows that
$g_d\inv\dots g_2\inv tg_2\dots g_d=t\in B$.

To study the remaining factor $g_d\inv\dots g_2\inv ug_2\dots g_d$, we
consider the subgroup
\[
U_{\al_1;\al_2,\dots,\al_d}:=\<U_{i_1\al_1+\dots+i_d\al_d}\mid
i_1,\dots,i_d\in\Z,\;i_1>0\>
\]
(where $U_\ga$ denotes the trivial group whenever $\ga$ is not a root).

Fix $2\le k\le d$. Let $\be=i_1\al_1+\dots+i_d\al_d\in\Phi$ for some
$i_1,\dots,i_d\in\Z$, $i_1>0$. For all $u_\be\in U_\be$ and all $u_{\al_k}\in
U_{\al_k}$, a well known commutation rule (see~\cite[Proposition~8.2.3]{Sp})
implies that
\[
u_{\al_k}u_\be u_{\al_k}\inv\in\<U_{i\be+j\al_k}\mid i>0,\,j\ge0\>\subset
U_{\al_1;\al_2,\dots,\al_d}.
\]
Similarly, $u_{-\al_k}u_\be u_{-\al_k}\inv\in U_{\al_1;\al_2,\dots,\al_d}$
for all $u_{-\al_k}\in U_{-\al_k}$. Using again that
$L_{\al_k}=\<U_{\al_k},U_{-\al_k}\>$, it follows that $g_k\inv
U_{\al_1;\al_2,\dots,\al_d}g_k=U_{\al_1;\al_2,\dots,\al_d}$.

In particular, $g_d\inv\dots g_2\inv ug_2\dots g_d\in
U_{\al_1;\al_2,\dots,\al_d}$. To complete the proof, it remains to show that
$U_{\al_1;\al_2,\dots,\al_d}\subset B$, or, in other words, that $\sum_k
i_k\al_k$ cannot be a negative root if $i_1>0$. Write $\|\al\|^2:=(\al|\al)$
for all $\al$ in the root lattice. By orthogonality, we have
\[
\|i_1\al_1+\dots+ i_d\al_d\|^2=i_1^2\,\|\al_1\|^2+\dots+i_d^2\,\|\al_d\|^2,
\]
so there are two cases:
\begin{itemize}
\item if $\al_1$ is long (in its component), then $\sum_ki_k\al_k$ cannot be
a root (except for $\al_1$ itself);
\item if $\al_1$ is short, then $\sum_ki_k\al_k$ can only be a root if it is
of the form $\al_1\pm\al_k$ for some $2\le k\le d$. But by assumption,
$\al_1\not<\al_k$, so $\al_1+\al_k$ and $\al_1-\al_k$ are either positive
roots (when $\al_1>\al_k$) or nonroots (when $\al_1,\al_k$ are incomparable).
\end{itemize}
\end{proof}
\begin{cor}
The set $E(C)$ is a subvariety of the flag variety $G/B$.
\end{cor}
\begin{proof}
Use Proposition~\ref{prodP1} and the fact that the $L_{\al_k}/B_{\al_k}$ are
complete.
\end{proof}
\begin{rem}
\label{indepnumb}
The subvariety $E(C)$ does not depend on the numbering of the elements of
$\Om$ (as long as one chooses a nonincreasing one).
\end{rem}
\begin{proof}
Let $\be_1,\dots,\be_d$ be another nonincreasing numbering of the elements of
$\Om$.
For any incomparable roots $\ga,\ga'$, the commutator $(L_\ga,L_{\ga'})$ is
trivial, so it is enough to show that the sequence $\al_1,\dots,\al_d$ can be
rearranged into $\be_1,\dots,\be_d$ by successively swapping adjacent pairs
of incomparable roots.

We have $\be_1=\al_k$ for some $1\le k\le d$, and by assumption, this element
must be maximal in $\Om$. Therefore, $\al_k$ is incomparable with each of
$\al_1,\dots,\al_{k-1}$, so we may move it past these roots to the beginning
of the sequence. Now we have got two sequences with a common first term;
discarding it, we may reapply the same procedure inductively.
\end{proof}


\section{The embedding property and $T$-stability}

We retain all previous notation and assumptions.
\begin{theo}
\label{maintheo}
The map $j$ of Proposition~\ref{prodP1} is an embedding, and its image
$E(C)$ is a $T$-stable subvariety in $G/B$, isomorphic to a
product of $d$ projective lines.
\end{theo}
\begin{proof}
First, note that if $E(C)$ is $T$-stable, then for every $w'\in W$,
$E(w'C)=w'E(C)$ is also $T$-stable; we may therefore restrict the proof to the
case where $\dot w=1$.

Since for each root $\al$, the quotient $L_\al/B_\al$ is a projective line,
the last statement follows from the first one and from
Proposition~\ref{prodP1}. Now choose a representative $\dot s_\al\in\Norm_G(T)$
for the reflection $s_\al$ and recall that the Bruhat decomposition
$L_\al=B_\al\cup U_\al s_\al B_\al=B_\al\cup s_\al U_{-\al} B_\al$ induces an
open covering
\[
L_\al/B_\al
=\{uB_\al\mid u\in U_{-\al}\}\cup\{\dot s_\al uB_\al\mid u\in U_{-\al}\}
\]
of $L_\al/B_\al$ by two affine lines.

Taking products, we get an open
covering of $\prod_k(L_{\al_k}/B_{\al_k})$ by $2^d$ affine sets. More
precisely, for each subset $K\subset\{1,\dots,d\}$, we define
\[
V_K:=\{(\dot r_{K,1}u_1B_{\al_1},\dots,\dot r_{K,d}u_dB_{\al_d})\mid u_k\in
U_{-\al_k}\,\forall k\}\subset\prod_k(L_{\al_k}/B_{\al_k}),
\]
where
\[
r_{K,k}:=
\begin{cases}
s_{\al_k}&\text{if $k\in K$,}
\\
1&\text{if $k\not\in K$.}
\end{cases}
\]
The set $V_\emptyset$ is dense in $\prod_k(L_{\al_k}/B_{\al_k})$, so its
image
\[
j(V_\emptyset)=\{u_1\dots u_dB\mid u_k\in U_{-\al_k}\,\forall k\}
\]
is also dense in
$j\Bigl(\prod_k(L_{\al_k}/B_{\al_k})\Bigr)=E(C)$. Since $T$
normalizes each root group $U_{-\al_k}$, it follows that
$E(C)$ is $T$-stable.

It remains to show that $j$ is an embedding. Since the
variety $\prod_k(L_{\al_k}/B_{\al_k})$ is complete, it is enough to show
that
\begin{itemize}
\item[(i)] the restriction of $j$ to each open affine set $V_K$ is an
embedding, and
\item[(ii)] $j$ is injective.
\end{itemize}

Condition (i) holds for $K=\emptyset$: indeed, the multiplication map
$U_{-\al_1}\times\dots\times U_{-\al_d}\to U^{-}$ is well known to be an
embedding (see e.g.~\cite[\S8.2.1]{Sp}), and so is the canonical map
$U^{-}\to G/B:u\mapsto uB$ (by the Bruhat decomposition of $G$).

To show condition (i) for an arbitrary subset $K\subset\{1,\dots,d\}$, we
factor the restriction of $j$ to $V_K$ through an isomorphism $V_K\simeq
V_\emptyset$, as follows. Consider an element
$(\dot r_{K,1}u_1B_{\al_1},\dots,\dot r_{K,d}u_dB_{\al_d})\in V_K$. Recall that for any
two orthogonal roots $\al,\be$, $s_\al$ normalizes $U_\be$ (see
e.g.~\cite[\S9.2.1]{Sp}), so we may rewrite
\[
\dot r_{K,1}u_1\dots \dot r_{K,d}u_d=(\dot r_{K,1}\dots \dot r_{K,d})(u'_1\dots
u'_d)=(\prod_{k\in K}\dot s_{\al_k})(u'_1\dots u'_d)
\]
for some $u'_1\in U_{-\al_1},\dots,u'_d\in U_{-\al_d}$. It is clear that the
map
\[
\sigma_K:V_K\to V_\emptyset:
(\dot r_{K,1}u_1B,\dots,\dot r_{K,d}u_dB)\mapsto(u'_1B,\dots,u'_dB)
\]
is an isomorphism, and by construction, the restriction of $j$ to $V_K$ is
equal to the composition $V_K\xto{\sigma_K}V_\emptyset\xto j G/B\to G/B$,
where the last arrow is (left) multiplication by $\prod_{k\in K}\dot
s_{\al_k}$. This shows condition~(i).

Finally, we must check condition~(ii), so let $p=(p_1,\dots,p_d)$ and
$q=(q_1,\dots,q_d)$ be two points in $\prod_k(L_{\al_k}/B_{\al_k})$. Let
$K\subset\{1,\dots,d\}$ be maximal (w.r.t.\ set inclusion) such that $p\in
V_K$: then for each $k\in K$, $p_k=\dot s_{\al_k}u_kB_{\al_k}$ with $u_k\in
U_{-\al_k}$, and for each $k\not\in K$, we must have $p_k=B_{\al_k}$. It
follows that
\[
j(p)=\left(\prod_{k\in K}\dot s_{\al_k}u_k\right)B=\left(\prod_{k\in
K}u^{+}_k\right)\left(\prod_{k\in K}\dot s_{\al_k}\right)B
\]
for some $u^{+}_1\in U_{\al_1}$, $\dots$, $u^{+}_d\in U_{\al_d}$. Therefore,
$j(p)$ lies in the Schubert cell $B\left(\prod_{k\in K}s_{\al_k}\right)B/B$.
Similarly, if $L\subset\{1,\dots,d\}$ is maximal such that $q\in V_L$, then
$j(q)\in B\left(\prod_{k\in L}s_{\al_k}\right)B/B$. Now assume $p\ne q$; there
are two cases:
\begin{itemize}
\item if $K\ne L$, then $j(p)\ne j(q)$ because they lie in
different Schubert cells of $G/B$;
\item if $K=L$, then $j(p)\ne j(q)$ by condition~(i).
\end{itemize}
This shows condition~(ii).
\end{proof}
\begin{rem}
\label{indepw}
The variety $E(C)$ does not depend on the choice of an element $w\in C$, nor of
its representative $\dot w\in\Norm_G(T)$.
\end{rem}
\begin{proof}
The last part follows from $T$-stability (see Theorem~\ref{maintheo}). For
the first part, recall again~\cite[\S9.2.1]{Sp} that if $\al,\be$ are
orthogonal roots, then $s_\al U_\be s_\al\inv=U_\be$; since
$L_\be=\<U_\be,U_{-\be}\>$, we then also have $s_\al L_\be s_\al\inv=L_\be$.
Now use this fact in the definition of $E(C)$.
\end{proof}

\section{Two examples}

\begin{ex}
$G=\SL(n)$ acts on the set of all (full) flags $F_1\subset
F_2\subset\dots\subset F_{n-1}$ (with $\dim F_i=i-1$) of linear subspaces in
the projective space $\P^{n-1}$. Choose a ``hypertetrahedron'' in $\P^{n-1}$,
i.e.\ $n$ linearly independent points $p_1,\dots,p_n\in\P^{n-1}$. Let $B$ be
the stabilizer of the flag
\[
\<p_1\>\subset\<p_1,p_2\>\subset\dots\subset\<p_1,\dots,p_{n-1}\>
\]
(where $\<\quad\>$ denotes linear span in $\P^{n-1}$) and let $T$ be the
simultaneous stabilizer of the vertices $p_1,\dots,p_n$. Identifying $W$ with
the symmetric group $S_n$, the $T$-fixed points in $G/B$ are the flags of the
form
\[
\<p_{w(1)}\>\subset\<p_{w(1)},p_{w(2)}\> \subset\dots\subset
\<p_{w(1)},\dots,p_{w(n-1)}\>
\]
for some $w\in W$.

Let $\{\al_1,\dots,\al_d\}$ be a set of positive and pairwise orthogonal
roots; each $s_{\al_k}$ is a transposition $(a_k\,b_k)\in S_n$, and the
sequence is automatically nonincreasing. For each $k$, pick a point $q_k$ on
the projective line $\<p_{a_k},p_{b_k}\>$; to the tuple $(q_1,\dots,q_d)$, we
associate the flag $F(q_1,\dots,q_d)$ whose $(i-1)$-dimensional component is
obtained from the expression $\<p_1,\dots,p_i\>$ by replacing $p_{a_k}$ by
$q_k$ whenever $a_k\le i<b_k$. For example, when $n=4$ and $d=2$, the flag
$F(q_1,q_2)$ is of one of the following forms:
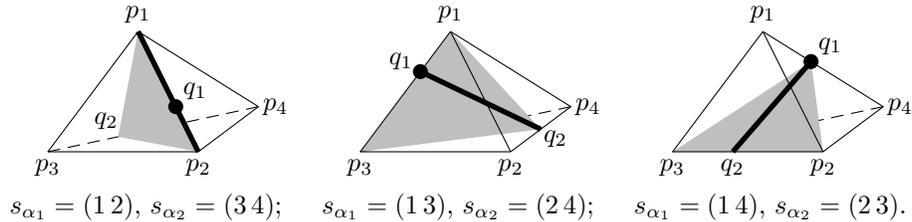
\begin{center}
\begin{tabular}{ccc}
\psset{unit=2mm,linewidth=0.4pt}
\begin{pspicture}(-3,-2)(16,11)
\qline(0,0)(10,0)
\qline(0,0)(6,8)
\qline(10,0)(6,8)
\qline(10,0)(14,3)
\qline(14,3)(6,8)
\psline[linestyle=dashed](0,0)(14,3)
\pspolygon*[linecolor=lightgray](6,8)(10,0)(4.667,1)
\psline[linewidth=2pt](10,0)(6,8)
\pscircle*(8.5,3){0.5}
\rput[b]{0}(6,8.6){$p_1$}
\rput[t]{0}(10,-0.5){$p_2$}
\rput[t]{0}(0,-0.5){$p_3$}
\rput[l]{0}(14.3,3){$p_4$}
\rput[bl]{0}(9,3.5){$q_1$}
\rput[br]{0}(4.667,1.5){$q_2$}
\end{pspicture}
&
\psset{unit=2mm,linewidth=0.4pt}
\begin{pspicture}(-3,-2)(16,11)
\qline(10,0)(14,3)
\qline(14,3)(6,8)
\psline[linestyle=dashed](0,0)(14,3)
\pspolygon*[linecolor=lightgray](6,8)(0,0)(12,1.5)
\psline[linewidth=2pt](4,5.333)(12,1.5)
\pscircle*(4,5.333){0.5}
\qline(0,0)(10,0)
\qline(0,0)(6,8)
\qline(10,0)(6,8)
\rput[b]{0}(6,8.6){$p_1$}
\rput[t]{0}(10,-0.5){$p_2$}
\rput[t]{0}(0,-0.5){$p_3$}
\rput[l]{0}(14.3,3){$p_4$}
\rput[br]{0}(3.5,5.5){$q_1$}
\rput[tl]{0}(12.2,1.4){$q_2$}
\end{pspicture}
&
\psset{unit=2mm,linewidth=0.4pt}
\begin{pspicture}(-3,-2)(16,11)
\qline(10,0)(14,3)
\qline(14,3)(6,8)
\psline[linestyle=dashed](0,0)(14,3)
\pspolygon*[linecolor=lightgray](0,0)(10,0)(9.2,6)
\psline[linewidth=2pt](9.2,6)(4,0)
\pscircle*(9.2,6){0.5}
\qline(0,0)(10,0)
\qline(0,0)(6,8)
\qline(10,0)(6,8)
\rput[b]{0}(6,8.6){$p_1$}
\rput[t]{0}(10,-0.5){$p_2$}
\rput[t]{0}(0,-0.5){$p_3$}
\rput[l]{0}(14.3,3){$p_4$}
\rput[bl]{0}(9.7,6.5){$q_1$}
\rput[t]{0}(4,-0.5){$q_2$}
\end{pspicture}
\\
$s_{\al_1}=(1\,2)$, $s_{\al_2}=(3\,4)$;
&
$s_{\al_1}=(1\,3)$, $s_{\al_2}=(2\,4)$;
&
$s_{\al_1}=(1\,4)$, $s_{\al_2}=(2\,3)$.
\end{tabular}
\end{center}
\vspace*{5mm}

Then $E(\<s_{\al_1},\dots,s_{\al_d}\>)$ is the set of all
$F(q_1,\dots,q_d)$, each $q_k$ varying on the line $\<p_{a_k},p_{b_k}\>$. If
$w\in W\simeq S_n$, then $E(w\<s_{\al_1},\dots,s_{\al_d}\>)$ is obtained from
the above description by replacing each $p_i$ by $p_{w(i)}$.

(Note that the embedding property of Theorem~\ref{maintheo} has an easy proof
here: indeed, the $s_{\al_k}$ pairwise commute, so the sets $\{a_k,b_k\}$
are pairwise disjoint, and therefore the projective lines
$\<p_{a_k},p_{b_k}\>$ on which the $q_k$ vary are pairwise skew.)
\end{ex}
\begin{ex}
$G=\Sp(4)$ acts on the set of all isotropic flags in $\P^3$, i.e.\ pairs
$(p,\ell)$ with $\ell\subset\P^3$ an isotropic line (w.r.t.\ the symplectic
form defining $\Sp(4)$) and $p$ a point on $\ell$. Choose a (skew)
``isotropic square'' in $\P^3$, i.e.\ four points $p_1,p_2,p_3,p_4\in\P^3$
such that all lines $\ell_{ij}:=\<p_i,p_j\>$ are isotropic except $\ell_{13}$
and $\ell_{24}$. Let $B$ be the stabilizer of the flag $(p_1,\ell_{12})$ and
$T$ be the simultaneous stabilizer of the vertices $p_1,p_2,p_3,p_4$. The
$T$-fixed points in $G/B$ are the flags of the form $(p_i,\ell_{ij})$ (with
$ij\ne13$ and $ij\ne24$).

Let $\al,\be$ be two positive orthogonal roots. There are two cases.
\begin{itemize}
\item If $\al,\be$ are both short, with $\al>\be$, then
$E(\<s_\al,s_\be\>)$ is the set of all isotropic flags $(p,\ell)$ such that
$\ell$ crosses $\ell_{14}$ and $\ell_{23}$.
\item If $\al,\be$ are both long (hence incomparable), then
$E(\<s_\al,s_\be\>)$ is the set of all isotropic flags $(p,\ell)$ such that
$p$ lies on the (nonisotropic) line $\ell_{13}$, and $\ell$ also crosses the
(nonisotropic) line $\ell_{24}$.
\end{itemize}
Again, if $w\in W$, then $E(w\<s_\al,s_\be\>)$ is obtained from this
description by applying $w$, viewed as a ``symmetry'' of
the ``square'' $p_1p_2p_3p_4$.
\end{ex}

\section{$T$-orbits in $E(C)$}

Let $C$ be an orthocell of rank $d$ and identify
$E(C)\simeq\P^1\times\dots\times\P^1$ via the map $j$ of
Proposition~\ref{prodP1}.
\begin{prop}
The $T$-action on $E(C)$ has $3^d$ orbits, viz.\ the subsets of the form
$A_1\times\dots\times A_d$, where each $A_k\subset\P^1$ is one of $\{0\}$,
$\{\infty\}$, $\P^1\setminus\{0,\infty\}$.
\end{prop}
\begin{proof}
For each $\al\in\Phi$, choose an isomorphism $u_\al:(\C,{+})\to U_\al$ such
that $tu_\al(z)t\inv=u_\al(\al(t)z)$ for all $t\in T$ and all
$z\in\C$~\cite[Proposition~8.1.1(i)]{Sp}.

Consider first the case where $C=\<s_{\al_1},\dots,s_{\al_d}\>$ (i.e.\ $w=1$).
The set of all $u_{-\al_1}(z_1)\dots u_{-\al_d}(z_d)B$,
$(z_1,\dots,z_d)\in\C^n$, is an open dense subset of $E(C)$ (cf.\ the proof of
Theorem~\ref{maintheo}), and we have
\[
t\,u_{-\al_1}(z_1)\dots u_{-\al_d}(z_d)B
=u_{-\al_1}(\al_1(t)\inv z_1)\dots u_{-\al_d}(\al_d(t)\inv z_d)B.
\]
By orthogonality, the $\al_k$ are linearly independent, hence the morphism
\[
T\to\C^*\times\dots\times\C^*:t\mapsto(\al_1(t),\dots,\al_d(t))
\]
is surjective. Therefore, by continuity, this morphism turns the $T$-action
on $E(C)$ into the natural componentwise action of
$\C^*\times\dots\times\C^*$ on $\P^1\times\dots\times\P^1$. The orbit
structure is now clear.

The general case $C=w\<s_{\al_1},\dots,s_{\al_d}\>$ is obtained similarly, by
multiplying $\dot wu_{-\al_1}(z_1)\dots u_{-\al_d}(z_d)B$ by $wtw\inv$.
\end{proof}
To describe the $T$-orbit closures, define a \emph{subcell} of an orthocell
$C=w\<s_{\al_1},\dots,s_{\al_d}\>$ to be an orthocell of the form
$w'\<s_{\al'_1},\dots,s_{\al'_e}\>$ for some
$\{\al'_1,\dots,\al'_e\}\subset\{\al_1,\dots,\al_d\}$ and some $w'\in C$.
\begin{cor}
The $T$-orbit closures in $E(C)$ are exactly the $E(C')$ with $C'$ a subcell of
$C$.
\end{cor}

\end{document}